\documentclass[10pt,reqno]{amsart}

\usepackage{amssymb}
\usepackage{xypic}
\xyoption{all}
\def\umono{\ar@{_{(}->}[u]}
\def\uumono{\ar@{_{(}->}[uu]}

\def\lmono{\ar@{_{(}->}[l]}
\def\llmono{\ar@{_{(}->}[ll]}

\newcommand{\Z}{{\mathbb Z}}

\newcommand{\F}{{\mathbb F}}

 \renewcommand{\P}[1]{\mathcal{P}^{#1}}

\newcommand{\map}{\operatorname{map}\nolimits}

\renewcommand{\ker}{\operatorname{Ker}\nolimits}
\newcommand{\im}{\operatorname{Im}\nolimits}
\newcommand{\coker}{\operatorname{Coker}\nolimits}
\newcommand{\A}{\ifmmode{\mathcal{A}}\else${\mathcal{A}}$\fi}
\newcommand{\K}{\ifmmode{\mathcal{K}}\else${\mathcal{K}}$\fi}
\newcommand{\U}{\ifmmode{\mathcal{U}}\else${\mathcal{U}}$\fi}
\newcommand{\T}{\ifmmode{\mathcal{T}}\else${\mathcal{T}}$\fi}
\newcommand{\FF}{\ifmmode{\mathcal{F}}\else${\mathcal{F}}$\fi}
\newcommand{\LL}{\ifmmode{\mathcal{L}}\else${\mathcal{L}}$\fi}

\newtheorem{theorem}{Theorem}[section]
\newtheorem{proposition}[theorem]{Proposition}
\newtheorem{corollary}[theorem]{Corollary}
\newtheorem{lemma}[theorem]{Lemma}
\newtheorem{definition}[theorem]{Definition}
\newtheorem{remark}[theorem]{Remark}
\newtheorem{example}[theorem]{Example}

\title[On the cohomology of highly connected covers of finite complexes]
{On the cohomology of highly connected covers of finite complexes}

\author{Nat\`{a}lia Castellana}
\author{Juan A. Crespo}
\author{J\'er\^{o}me Scherer}

\thanks{All three authors are partially supported by MEC grant MTM2004-06686.
The first and third authors were partially supported by the
Mittag-Leffler Institute in Sweden. The third author is supported
by the program Ram\'on y Cajal, MEC, Spain.}


\begin{document}


\begin{abstract}
Relying on the computation of the Andr\'e-Quillen homology groups
for unstable Hopf algebras, we prove that the mod $p$ cohomology
of the $n$-connected cover of a finite $H$-space is always
finitely generated as algebra over the Steenrod algebra.
\end{abstract}


\maketitle


\section*{Introduction}
\label{sec intro}

Consider the $n$-connected cover of a finite complex. Does its
(mod $p$) cohomology satisfy some finiteness property? Such a
question has already been raised by McGibbon and M\o ller in
\cite{MR1444714}, but no satisfactory answer has been proposed. We
do not ask here for an algorithm which would allow to make
explicit computations. We rather look for a general
structural statement which would tell us to what kind of class
such cohomologies belong. The prototypical theorems we have in
mind are the Evens-Venkov result, \cite{MR0137742},
\cite{MR0108788}, that the cohomology of a finite group is
Noetherian, the analog for $p$-compact groups obtained by Dwyer
and Wilkerson \cite{MR1274096}, and the fact that the mod $p$
cohomology of an Eilenberg-Mac Lane space $K(A, n)$, with $A$
abelian of finite type, is finitely generated as an algebra over
the Steenrod algebra, which can easily been inferred from the work
of Serre \cite{MR0060234} and Cartan \cite{MR0087934}.

This last observation leads us to ask first whether or not the mod
$p$ cohomology of a finite Postnikov piece is also finitely
generated as an algebra over the Steenrod algebra and second,
since a finite complex $X$ and its $n$-connected cover $X \langle
n \rangle$ only differ in a finite number of homotopy groups, if
$H^*(X \langle n \rangle; \F_p)$ satisfies the same property. We
offer in this paper a positive answer when $X$ is an $H$-space,
based on the analysis of the fibration $P \rightarrow X \langle n
\rangle \rightarrow X$, where $P$ is a finite Postnikov piece. In
fact we prove a strong closure property for $H$-fibrations.

\medskip

{\bf Theorem~\ref{generalclosure}.}
\noindent {\it Let $F \rightarrow E \rightarrow B$ be an
$H$-fibration in which both $H^*(F; \F_p)$ and $H^*(B; \F_p)$ are
finitely generated unstable algebras. Then so is $H^*(E; \F_p)$.}

\medskip

This applies in particular to highly connected covers of finite
$H$-spaces, see Theorem~\ref{covers}. In our previous work
\cite{deconstructing} we proved that the theorem holds whenever
the base space is an Eilenberg-Mac Lane space. The proof relied
mainly on Smith's work \cite{MR0275435} on the Eilenberg-Moore
spectral sequence.

Our starting point here is the same and we need to analyze
carefully certain Hopf subalgebras of $H^*(F; \F_p)$. Observe that
the property for an unstable algebra $K$ to be a finitely
generated $\A_p$-algebra is equivalent to say that the module of
the undecomposable elements $QK$ is finitely generated as unstable
module. It is often more handy to work with this module because it
is smaller than the whole algebra and, above all, the category of
unstable modules is locally Noetherian, \cite{MR868302}.

The main problem (or interest) with the functor $Q(-)$ is the
failure of left exactness. To what extent this functor is not left
exact is precisely measured by Andr\'e-Quillen homology
$H^Q_*(-)$. In our setting we keep control of the size of these
unstable modules.

\medskip

{\bf Proposition~\ref{H1isfg}.}
\noindent {\it Let $A$ be a Hopf algebra which is a finitely
generated unstable ${\mathcal A}_p$-algebra. Then $H_0^Q(A) = QA$
and $H^Q_1(A)$ are both finitely generated unstable modules.}

\medskip

As the higher groups are all trivial (see Proposition \ref{H1asvsprop}), this gives a quite accurate
description of Andr\'e-Quillen homology in our situation. The
relevance of Andr\'e-Quillen homology in homotopy theory is
notorious since Miller solved the Sullivan conjecture,
\cite{Miller}. The typical result which is needed in his work, and
which has been then extended by Lannes and Schwartz,
\cite{MR827370}, is that the module of indecomposable elements of
an unstable algebra $K$ is locally finite if and only if so are
\emph{all} Andr\'e-Quillen homology groups of $K$.

Proposition~\ref{H1isfg} yields then our main algebraic structural
result about the category of unstable Hopf algebras.

\medskip

{\bf Theorem~\ref{locallynoetherian}.}
\noindent {\it Let $B$ be a Hopf algebra which is a finitely
generated unstable ${\mathcal A}_p$-algebra. Then so is any
unstable Hopf subalgebra.}

\medskip

\noindent{\bf Acknowledgments.} We would like to thank Jean Lannes
for his guided reading of \cite{MR1156389} at the Mittag-Leffler
Institute and Carles Broto for many discussions on this material.

\section{Andr\'e-Quillen homology of Hopf algebras}
\label{section AndreQuillen}
In this section, we compute Andr\'e-Quillen homology for Hopf
algebras, and introduce the action of the Steenrod algebra in the
next one. A clear and short introduction to Andr\'e-Quillen
homology can be found in Bousfield's \cite[Appendix]{MR0380779},
see also Goerss' book \cite{MR1089001}.

Let us briefly recall from Schwartz's book \cite{MR95d:55017} how
one computes Andr\'e-Quillen homology in our setting. The
symmetric algebra comonad $S(-)$ yields a simplicial resolution
$S^\bullet(A)$ for any commutative algebra $A$. The
Andr\'e-Quillen homology group $H^Q_i(A)$ is the $i$-th homology
group of the complex obtained from $S^\bullet(A)$ by taking the
module of indecomposable elements (and the differential is the
usual alternating sum). This is a graded $\F_p$-vector space. Long
exact sequences arise from certain extensions, just like in the
dual situation for the primitive functor,
\cite[Theorem~3.6]{MR0380779}.

\begin{lemma}
\label{longexact}
Let $A$ be a Hopf subalgebra of a Hopf algebra $B$ of finite type.
Then there is a long exact sequence
$$
\dots \rightarrow H_2^Q(B//A) \rightarrow  H_1^Q(A) \rightarrow
H_1^Q(B) \rightarrow H_1^Q(B//A) \rightarrow QA \rightarrow QB
\rightarrow Q(B//A)
$$
in Andr\'e-Quillen homology.
\end{lemma}

\begin{proof}
Long exact sequences in Andr\'e-Quillen homology are induced by
cofibrations of simplicial algebras. However, the inclusion $A
\subset B$ of a sub-Hopf algebra is not a cofibration in general
(seen as a constant simplicial object). To get around this
difficulty we use Goerss' argument from
\cite[Section~10]{MR1089001} and we reproduce it here in our
framework.

For any morphism $f: A \rightarrow B$ of simplicial algebras,
there is a spectral sequence, \cite[Proposition~4.7]{MR1089001},
$\hbox{\rm Tor}^{\pi_*A}_p(\F_p, \pi_* B)_q$ converging to the
homotopy groups $\pi_{p+q} Cof(f)$ of the homotopy cofiber. Now,
because $B$ is of finite type, it is always a free $A$-module by
the Milnor-Moore result \cite[Theorem~4.4]{MR0174052}. Thus the
$E_2$-term is isomorphic to $\hbox{\rm Tor}^{\F_p}_0(\F_p, B//A)_*
\cong B//A$. The spectral sequence collapses and hence $Cof(f)$ is
weakly equivalent to $B//A$. In particular, we have the desired
long exact sequence.
\end{proof}

Following the terminology used in \cite[Section~6]{MR0275435}, we
introduce the following definition.

\begin{definition}
\label{defcoexact}
{\rm A sequence of (Hopf) algebras
$$
\F_p \longrightarrow A \longrightarrow B \longrightarrow C
\longrightarrow \F_p
$$
is \emph{coexact} if the morphism $A \rightarrow B$ is a
monomorphism and its cokernel $B//A$ is isomorphic to $C$ as a
(Hopf) algebra.}
\end{definition}

We can thus restate the previous lemma by saying that coexact
sequences of Hopf algebras induce long exact sequences in
Andr\'e-Quillen homology.

By the Borel-Hopf decomposition theorem
\cite[Theorem~7.11]{MR0174052}, any Hopf algebra of finite type is
isomorphic, as an algebra, to a tensor product of monogenic Hopf
algebras, i.e. either a truncated polynomial algebra of the form
$\F_p[x_i]/(x_i^{p^{k_i}})$, where $p^{k_i}$ is the \emph{height}
of the generator $x_i$, or a polynomial algebra of the form
$\F_p[y_j]$, or, when $p$ is odd, an exterior algebra
$\Lambda(z_i)$. Let us denote by $\xi$ the Frobenius map, sending
an element $x$ to its $p$-th power $x^p$.

\begin{proposition}
\label{H1asvsprop}
Let $A$ be a Hopf algebra of finite type. Then $H_0^Q(A) = QA$ and
$H^Q_1(A)$ is isomorphic to the $\F_p$-vector space generated by
the elements $\xi^{k_i} x_i$ of degree $p^{k_i} \cdot |x_i|$ where $x_i\in A$ is a generator of height $p^{k_i}$, $0<k_i<\infty$ .
Moreover $H^Q_n(A) = 0$ if $n \geq 2$.
\end{proposition}

\begin{proof}
Consider the symmetric algebra $S(QA)$ and construct an algebra
map $S(QA) \twoheadrightarrow A$ by choosing representatives in
$A$ of the indecomposable elements. We have then a coexact
sequence of algebras $\F_p[\xi^{k_i}x_i] \hookrightarrow S(QA)
\twoheadrightarrow A$ and $A$, as Hopf algebra, can be seen as the
quotient $S(QA)// \F_p[\xi^{k_i}x_i]$. Since $S(QA)$ is a free
commutative algebra, $H_n^Q\bigl( S(QA) \bigr) =0$ for all $n \geq
1$. Likewise $H_n^Q\bigl( \F_p[\xi^{k_i}x_i] \bigr) =0$ for all $n
\geq 1$. Now, Lemma~\ref{longexact} allows us to identify
$H_1^Q(A) \cong H_0^Q(\F_p[\xi^{k_i}x_i]) \cong \oplus_k \F_p
\langle \xi^{k_i} x_i \rangle$, as a graded vector space.
\end{proof}

The vanishing of the higher Andr\'e-Quillen homology groups, or in
other words the fact that the functor $Q(-)$ has homological
dimension $\leq 1$ for Hopf algebras, has been analyzed by
Bousfield in the dual situation \cite[Theorem~4.1]{MR0380779}.
The next lemma is now a reformulation of the preceding
proposition.

\begin{lemma}
\label{H1asvs}
Let $A$ be a Hopf algebra of finite type and denote by $x_i$ the
truncated polynomial generators. Then $H^Q_1(A)$ is isomorphic to
the $\F_p$-vector space generated by the elements $x_i^{\otimes
p^{k_i}}\in S(A)$, where $p^{k_i}$ is the height of~$x_i$.
\end{lemma}

\begin{proof}
We have to compute the first homology group of the complex
$$
\dots
\rightarrow S^2(A) \xrightarrow{d} S(A) \xrightarrow{m} A \,.
$$
The morphisms are given by the alternating sums of the face maps.
Let us use the symbols $\otimes$ for the tensor product in $S(A)$
and $\boxtimes$ for the next level in $S(S(A))$. If $\eta_A: S(A)
\rightarrow A$ is the counit defined by $\eta_A(a \otimes b) =
ab$, the two face maps $S^2(A) \rightarrow S^1(A)$ are then
$S(\eta_A)$ and $\eta_{S(A)}$.

Thus $m(a) = S(\eta_A)(a) - \eta_{S(A)} = a-a = 0$ and $m(a
\otimes b) = ab$, since $\eta_{S(A)}(a \otimes b) = a \otimes b$
is decomposable. Likewise $d(w) = \eta_A(w)$ on elements $w \in
S(A)$ and $d(v \boxtimes w) =  v \otimes w - \eta_A(v) \otimes
\eta_A(w)$ for $v, w \in S(A)$. The elements $x_i^{\otimes
p^{k_i}}$ clearly belong to the kernel of $m$. To compare them to
the generators $\{ \xi^{k_i} x_i \}$ of $H_0^Q(\F_p[\xi^{k_i}
x_i]) \cong H_1^Q(A)$, we apply $S^\bullet$ to the coexact
sequence of algebras $\F_p[\xi^{k_i}x_i] \hookrightarrow S(QA)
\twoheadrightarrow A$. The snake Lemma yields a connecting
morphism $\ker(m) \rightarrow H_0^Q(\F_p[\xi^{k_i}x_i])$, which
sends precisely $x_i^{\otimes p^{k_i}}$ to $\xi^{k_i} x_i$.
\end{proof}

\begin{remark}
\label{Goerss}
{\rm Alternatively, one could use the identification of the first
Andr\'e-Quillen homology group $H_1^Q(A)$ with the indecomposable
elements of degree~$2$ in $\hbox{\rm Tor}_A(\F_p,\F_p)$,
\cite[Section~10]{MR1089001}. As an $\F_p$-vector space it is
generated by the transpotence elements $[x_i^{p^{k_i}-1}|x_i]$.
When $p=2$ and $k_i = 0$, this is not technically speaking a
transpotence element, but still it is $[x_i|x_i]$ which appears in
degree $2$, see for example \cite[Section~29-2]{kane}.}
\end{remark}

\section{Bringing in the action of the Steenrod algebra}
\label{section unstable}

The results of the previous section apply to Hopf algebras which
are finitely generated as algebras over the Steenrod algebra: they
are of finite type. Our aim in this section is to identify the
action of the Steenrod algebra on the unstable module $H^Q_1(A)$.
Good references on Andr\'e-Quillen homology for unstable algebras
are \cite[Chapter~7]{MR95d:55017} and of course \cite{Miller},
which showed the importance of Andr\'e-Quillen homology in a
topological context.

The $\F_p$-vector space $H^Q_1(A)$ is equipped with an action of
$\A_p$ because the Steenrod algebra acts on the symmetric algebra
via the Cartan formula. This yields the same unstable module
$H^Q_1(A)$ as the derived functor computed with a resolution in
the category of unstable algebras, \cite{MR827370} and
\cite[Proposition~7.2.2]{MR95d:55017}.

As expected with this type of questions, the case when $p=2$ is
slightly simpler than the case when $p$ is odd. To write a unified
proof, we use the well-known trick \cite{MR868302} to consider, in
the odd-primary case, the subalgebra of $\A_p$ concentrated in
even degrees. If $M$ is a module over $\A_p$, the module $M'$
concentrated in even degree is defined by $(M')^{2n} = M^{2n}$ and
$(M')^{2n+1} = 0$. This is not an $\A_p$-submodule of $M$, but it
is a module over $\A_p$ on which the Bockstein $\beta$ acts
trivially. Hence it can be seen as a module over the algebra
$\A_p'$, the subalgebra of $\A_p$ generated by the operations
$\P{i}$. When $p=2$ we adopt the convention that $\A_2' = \A_2$,
$\U' = \U$, and write $\P{i}$ for~$Sq^i$. Like in
\cite[1.2]{MR95d:55017}, for a sequence $I=(\varepsilon_0, i_1,
\varepsilon_1, \dots, i_n, \varepsilon_n)$ where the
$\varepsilon_k$'s are 0 or 1, we write $\P{I}$ for the operation
$\beta^{\varepsilon_0} \P{i_1} \beta^{\varepsilon_1} \dots
\P{i_n}\beta^{\varepsilon_n}$.

 In \cite[Appendice~B]{MR868302}, Lannes and Zarati prove that the
category $\U$ is locally noetherian, which they do by reducing the
proof to the case of $\U'$. We use their computations in the following lemma, in fact the
explicit version from \cite[1.8]{MR95d:55017}.

\begin{lemma}
\label{evenpart}
Let $M$ be an unstable module which is finitely generated over
$\A_p$. Then so is the module $M'$, over $\A_p'$.
\end{lemma}

\begin{proof}
The statement is a tautology when $p=2$. Let us assume $p$ is an
odd prime. In the category $\U'$ of unstable modules concentrated
in even degrees, $F'(2n)$ is the free object on one generator
$\iota_{2n}$ in degree $2n$. We must show that $M'$ is a quotient
of a finite direct sum of such modules. As we know that $M$ is a
quotient of a finite direct sum of $F(n)$'s, it is enough to prove
the lemma for a free module $F(n)$.

A basis over $\F_p$ for the module $F(n)$ is given by the elements
$\P{I}\iota_n$ where $I$ is admissible with excess $e(I) \leq n$.
Define $F(n)'_k$ to be the span over $\F_p$ of the elements
$\P{I}\iota_n$ with $e(I)\geq k$. Then $F(n)'_k/F(n)'_{k+1}$ is
zero when $k+n$ is odd and it is generated by the images of the
elements $\P{I}\iota_n$ where the $k$ Bocksteins appear in the
first $k+1$ possible slots. In particular, $F(n)'$ is generated by
these elements as an $\A_p'$-module.
\end{proof}

The generators for $H_1^Q(A)$ will be related to certain elements
in $QA$ we describe next.

\begin{lemma}
\label{fgu'}
Let $A$ be a Hopf algebra which is a finitely generated unstable
${\mathcal A}_p$-algebra and let $N$ be the submodule of $QA$
generated by the truncated polynomial generators $x_i$. Then $N'$
is finitely generated in $\U'$. There exists an integer $d$ and a
finite set $\{ x_{k,i} \, | \, 1 \leq k \leq d, 1 \leq i \leq n_k
\}$ of generators such that any element in $N$ of height $p^k$
can be written $\sum_{i} \theta_{k,i} x_{k,i}$ for some
(admissible) operations $\theta_{k,i} \in \A_p'$.
\end{lemma}

\begin{proof}
Since $\U$ is locally noetherian,
\cite[Theorem~1.8.1]{MR95d:55017}, the unstable module $N$ is
finitely generated, being a submodule of $QA$.
Thus, by Lemma~\ref{evenpart}, $N'$ is finitely generated over
$\A_p'$. This implies in particular that the height of the
truncated generators is bounded by some integer $p^d$ (the action
of the Steenrod algebra on $x_i$ can only lower the height by the
formulas \cite[1.7.1]{MR95d:55017}).

For $1 \leq k \leq d$, write $N'(k)$ for the submodule of $N'$
generated by the $x_i$'s of height $p^k$ and choose
representatives in $A$ of generators $x_{k,i}$ for the module
$N'(k)$, with $1 \leq i \leq n_k$. Hence the finite set $\{
x_{k,i} \, | \, 1 \leq k \leq d, 1 \leq i \leq n_k \}$ generates
$N'$.
\end{proof}

The relation $x = \sum_{i} \theta_{k,i} x_{k,i}$ holds in the
module of indecomposable elements (in fact in~$N'$). Beware that
the same relation holds also in the algebra $A$, but only up to
decomposable elements.

\begin{proposition}
\label{H1isfg}
Let $A$ be a Hopf algebra which is a finitely generated unstable
${\mathcal A}_p$-algebra. Then $H_0^Q(A) = QA$ and $H^Q_1(A)$ are
both finitely generated unstable modules.
\end{proposition}

\begin{proof}
Lemma~\ref{H1asvs} allows us to identify $H_1^Q(A) \cong \oplus_k
\F_p \langle x_i^{\otimes p^{k_i}} \rangle$, as a graded vector
space. We must now identify the action of the Steenrod algebra.

We claim that the finite set of elements $x_{k,i}^{\otimes
p^{k_i}}$ generates $H_1^Q(A)$ as unstable module. More precisely
we show that the relation $x = \sum_{i} \theta_{k,i} x_{k,i}$ in
$QA$ yields a relation for $x^{\otimes p^{k_i}}$ in $H_1^Q(A)$. To
simplify the notation, let us assume that the height of $x$ is
$p^k$ and that the relation is of the form $x = \sum_{j}
\theta_{j} x_{j}$ for generators $x_j$ of the same height. The
relation for $x$ holds in $A$ up to decomposable elements which
must have lower height. But if $a^{p^k} = 0 = b^{p^k}$,
then
$$
d[ a^{\otimes p^k} \boxtimes b^{\otimes p^k} - (a\otimes
b)^{\boxtimes p^k}] = a^{\otimes p^k} \otimes b^{\otimes p^k} -
a^{p^k} \otimes b^{p^k} - (a \otimes b)^{\otimes p^k} +
(ab)^{\otimes p^k} = (ab)^{\otimes p^k}
$$
and hence the decomposable elements disappear in $H_1^Q(A)$.
Therefore $x^{\otimes p^{k}} = \bigl( \sum_{i} \theta_{j} x_{j}
\bigr)^{\otimes p^{k}}$ in $H_1^Q(A)$. The operations $\theta_{j}$
live in $\A_p'$, so that we have basically to perform the
following computation in the symmetric algebra: $(\P{n}
x)^{\otimes p^k} = \P{p^k n}(x^{\otimes p^k})$. There exist thus
operations $\Theta_j \in \A_p'$ such that
$$
x^{\otimes p^{k}} = \sum_j (\theta_j x_j)^{\otimes p^k} = \sum
\Theta_j (x_j^{\otimes p^k})
$$
and the claim is proven.
\end{proof}

\begin{theorem}
\label{locallynoetherian}
Let $B$ be a Hopf algebra which is a finitely generated unstable
${\mathcal A}_p$-algebra. Then so is any unstable Hopf subalgebra.
\end{theorem}

\begin{proof}
Consider an unstable Hopf subalgebra $A \subset B$ and the
quotient $B//A$. By Lemma~\ref{longexact}, we have an associated
exact sequence in Andr\'e-Quillen homology
$$
H_1^Q(B//A) \rightarrow QA \rightarrow QB,
$$
in which the unstable modules $QB$ and $H_1^Q(B//A)$ are finitely
generated by Proposition~\ref{H1isfg}. Thus so is $QA$.
\end{proof}

\begin{example}
\label{oddprimes}
{\rm Let us consider the Hopf algebra $B = H^*(K(\Z/p, 2))$. When
$p$ is odd it is the tensor product of a polynomial algebra
$\F_p[\iota_2, \beta \P{1} \beta \iota_2, \beta \P{p} \P{1} \beta
\iota_2, \cdots]$, concentrated in even degrees, with an exterior
algebra $\Lambda(\beta \iota_2, \P{1} \beta \iota_2, ...)$.

We consider the Hopf subalgebra $A$ given by the image of the
Frobenius $\xi$. This is the polynomial
subalgebra
$$
\F_p[(\iota_2)^p, (\beta \P{1} \beta \iota_2)^p, (\beta \P{p}\P{1}
\beta \iota_2)^p, \dots]
$$
The quotient $B//A$ has an exterior part and a truncated
polynomial part where all generators have height $p$. The module
of indecomposable elements $Q(B//A)$ is isomorphic to $QB$. It is
a quotient of $F(2)$, and thus generated, as an unstable module,
by a single generator $\iota_2$ in degree $2$. The submodule
concentrated in even degree is a module over $\A_p'$. It is
finitely generated as well, by Lemma~\ref{evenpart}, but one needs
two generators $\iota_2$ and $\beta \P{1} \beta \iota_2$. Explicit
computations of the action of the Steenrod algebra can be found in
\cite{MR2002g:55016}.

Therefore $H_1^Q(B//A)$ is an unstable module, which is generated
by the elements $\iota_2^{\otimes p}$ and $(\beta \P{1} \beta
\iota_2)^{\otimes p}$.}
\end{example}

\begin{remark}
\label{Henn}
{\rm For plain unstable algebras, Theorem~\ref{locallynoetherian}
is false, as pointed out to us by Hans-Werner Henn. Consider
indeed the unstable algebra
$$
H^*(\mathbf CP^\infty \times S^2; \mathbf F_p) \cong \mathbf
F_p[x] \otimes E(y)
$$
where both $x$ and $y$ have degree~$2$. Take the ideal generated
by $y$, and add $1$ to turn it into an unstable subalgebra. Since
$y^2=0$, this is isomorphic, as an unstable algebra, to $\mathbf
F_p \oplus \Sigma^2 \mathbf F_p \oplus \Sigma^2
\widetilde{H}^*(\mathbf CP^\infty; \mathbf F_p)$, which is not
finitely generated.}
\end{remark}

\section{$H$-fibrations over Eilenberg-Mac Lane spaces}
\label{section EML}
In the second part of this paper we turn now our attention to
topological applications of the Andr\'e-Quillen homology
computation we have done previously. More precisely, we concentrate on
$H$-fibrations.

\begin{definition}
\label{fgclosure}
{\rm An $H$-space $B$ satisfies the \emph{(strong) fg closure
property} if, for any $H$-fibration $F \rightarrow E \rightarrow
B$, the cohomology  $H^*(E)$ is a finitely generated unstable
algebra if (and only if) so is $H^*(F)$.}
\end{definition}

We prove in this section that Eilenberg-Mac Lane spaces enjoy the
strong fg closure property. In \cite{deconstructing} we
established the fg closure property for $K(A,n)$ with $n\geq 2$,
which was sufficient to our purposes there.

Given $n\geq 2$, consider a non-trivial $H$-fibration  $F
\xrightarrow{i} E \xrightarrow{\pi} K(A, n)$ where $A$ is either
$\Z/p$ or a Pr\"ufer group $\Z_{p^\infty}$. This situation has been extensively and carefully studied by L. Smith in \cite{MR0275435}.

The $E_2$-term of the Eilenberg-Moore spectral sequence is given
by $\hbox{\rm Tor}_{H^*(K(A,n))}(H^*(E),\F_p)$ and converges to
$H^*(F)$. Since we deal with an $H$-fibration,
\cite[Theorem~2.4]{MR0275435} applies and $E_2 \cong
H^*(E)//\pi^*\otimes \hbox{\rm Tor}_{H^*(K(A,n))
\backslash\backslash \pi^*}(\F_p,\F_p)$ as algebras, where
$H^*(K(A,n)) \backslash\backslash \pi^*$ is the Hopf subalgebra
kernel of $\pi^*$. The first differential  is $d_{p-1}$
\cite[Theorem~4.7]{MR0275435}. The same argument as in
\cite[Section~5]{MR0275435} (done for stable Postnikov pieces) is
valid in our situation as well, and has been in fact already used
in this setting, see \cite[Proposition~7.3$^*$]{MR0275435}: on
algebra generators the next differentials must be zero, so that
the spectral sequence collapses at $E_p$. This term is generated
by $H^*(E)//\pi^*= E_p^{0,*}$, $s^{-1,0}\coker \beta {\mathcal
P}_0 \subset E_p^{-1,*}$, and $s^{-1,0}QH^*(K(A,n))\backslash
\backslash \pi^{odd} \subset E_p^{-1,*}$, where $\beta {\mathcal
P}_0 \colon QH^{odd}(K(A,n))\rightarrow QH^{even}(K(A,n))$ is
defined by $\beta {\mathcal P}_0(x)=\beta \P{t}(x)$ with
$2t+1=|x|$.

The algebra structure is described in  \cite[Proposition
7.3$^*$]{MR0275435} by means of coexact sequences, see
Definition~\ref{defcoexact}.



\begin{proposition}\cite{MR0275435}
\label{coexact}
Let $n \geq 2$ and consider an $H$-fibration $F \xrightarrow{i} E
\xrightarrow{\pi} K(A, n)$ where $A$ is either $\Z/p$ or a
Pr\"ufer group $\Z_{p^\infty}$. Then there is a coexact sequence
of Hopf algebras
$$
\F_p \longrightarrow H^*(E)//\pi^* \xrightarrow{i^*} H^*(F)
\longrightarrow R \longrightarrow \F_p,
$$
and $R$ is described in turn by a coexact sequence of Hopf algebras
$$
\F_p \longrightarrow \Lambda \longrightarrow R \longrightarrow S
\longrightarrow \F_p,
$$
where $\Lambda$ is an exterior algebra generated by
$s^{-1,0}\coker \beta {\mathcal P}_0$, and $S \subseteq
H^*(K(A,n-1))$ is a Hopf subalgebra. \hfill{\qed}
\end{proposition}

\begin{example}
\label{sphere}
{\rm Let us see how the well-known cohomology of $S^3\langle 3
\rangle$ can be identified with these tools. Consider the
fibration
$$
S^3\langle 3
\rangle\stackrel{i}{\rightarrow} S^3\stackrel{\pi}{\rightarrow}
K(\Z,3)
$$
In this situation $H^*(S^3)//\pi^*=0$ and $H^*(S^3\langle 3
\rangle) \cong \Lambda \otimes S$ by Proposition~\ref{coexact}.
Recall that $H^*(K(\Z,3))=\F_p[\beta {\mathcal P}_k u_3: k\geq
1]\otimes E({\mathcal P}_k u_3: k\geq 0)$ where ${\mathcal
P}_k=\P{p^{k-1}}\cdots \P{1}$ and $\mathcal P_0 u_3 = u_3$. Here
the Hopf algebra kernel $H^*(K(\Z,3))\backslash \backslash \pi^*$
is $\F_p[\beta {\mathcal P}_k u_3: k\geq 1]\otimes E({\mathcal
P}_k u_3: k\geq 1)$. One sees next that the cokernel of $\beta
{\mathcal P}_p$ is $\{\beta \P{1} u_3\}$ and $S\subset
H^*(K(\Z,2)) $ is generated by $\P{1}\iota_2=\iota_2^p$. This
implies that $H^*(S^3\langle 3 \rangle) \cong \F_p[x_{2p}]\otimes
E(\beta x_{2p})$.}
\end{example}

Proposition~\ref{coexact} will allow us to improve
\cite[Theorem~6.1]{deconstructing}. We rely on the following
obvious lemma, which we will use again in the next section.

\begin{lemma}
\label{extension}
Consider an $H$-space $B$ and assume that there exists an
$H$-fibration $B' \rightarrow B \rightarrow B''$ such that both
$B'$ and $B''$ satisfy the (strong) fg closure property. Then so
does $B$.
\end{lemma}

\begin{proof}
Consider an $H$-fibration $F \rightarrow E \rightarrow B$ and
construct the pull-back diagram of fibrations
\[
\diagram E' \rto \dto_{p'} & E \dto^{p} \rto & B'' \ar@{=}[d] \cr
B' \rto & B \rto & B''
\enddiagram
\]
The homotopy fiber of $p'$ is $F$, which allows to conclude.
\end{proof}

\begin{theorem}
\label{complement}
Let $A$ be a finite direct sum of copies of cyclic groups $\Z/p^r$
and Pr\"ufer groups $\Z_{p^\infty}$, and $n \geq 2$. Consider an
$H$-fibration $F \xrightarrow{i} E \xrightarrow{\pi} K(A, n)$.
Then $H^*(F)$ is a finitely generated $\A_p$-algebra if and only
if so is $H^*(E)$.
\end{theorem}

\begin{proof}
If we consider the fibration of Eilenberg-Mac Lane spaces induced
by a group extension $A'\rightarrow A\rightarrow A$, we see from
Lemma~\ref{extension} that we can assume that $A=\Z/p$ or
$\Z_{p^\infty}$.

Since $H^*(K(A,n))$ is finitely generated as algebra over
$\mathcal A_p$, so is its image $\im(\pi^*)\subseteq H^*(E)$.
Hence, to prove the theorem, it is enough to show that the module
of indecomposable elements $Q(H^*(E)//\pi^*)$ is a finitely
generated $\A_p$-module if and only if so is $QH^*(F)$.

Let us now apply Lemma~\ref{longexact} to the coexact sequences
from Proposition~\ref{coexact}. The unstable Hopf algebra $S$ is
an unstable Hopf subalgebra of $H^*(K(A, n))$. Thus
Theorem~\ref{locallynoetherian} implies that $S$ is finitely
generated over $\A_p$. Moreover, the exterior algebra $\Lambda$ is
identified with $E(s^{-1,0} \coker{\beta \P{}_0})$, where the
cohomological operation $\beta \P{}_0$ is to be understood as an
operation from the odd degree part of $QH^*(K(A, n))$ to the even
degree part. The latter module is finitely generated by
Lemma~\ref{evenpart}. Hence the cokernel is finitely generated as
well, as a module over the Steenrod algebra. The exact sequence in
Andr\'e-Quillen homology for the coexact sequence involving $R$
and Proposition~\ref{H1isfg} show that both $QR$ and $H_1^Q(R)$
are finitely generated unstable modules. Finally, since $\mathcal
U$ is a locally noetherian category,
\cite[Theorem~1.8.1]{MR95d:55017}, the exact sequence
$$
H_1^Q(R) \rightarrow Q(H^*(E)//\pi^*) \rightarrow QH^*(F)
\rightarrow QR
$$
implies that $QH^*F$ is a finitely generated $\A_p$-module if and
only if so is $Q(H^*(E)//\pi^*)$.
\end{proof}

\begin{remark}
\label{Lannes}
{\rm Another approach to Theorem \ref{complement} is to dualize
the work of Goerss, Lannes, and Morel in
\cite[Section~2]{MR1156389}. They analyze the homology sequence
$$H_*(\Omega^2 K(A, n)) \rightarrow H_*(\Omega F) \rightarrow
H_*(\Omega E)\rightarrow H_*(\Omega K(A, n))$$ for a fibration of
spaces $F \rightarrow E \rightarrow K(A, n)$ and measure its
failure to be exact. This can be dualized and actually works for
$H$-fibrations, not only loop fibrations. Let us quickly sketch
the key ideas. Consider now an $H$-fibration $F \rightarrow E
\rightarrow K(A, n)$ and the complex
$$
H^*(K(A,n))\xrightarrow{\pi^*} H^* E \xrightarrow{i^*} H^* F
\longrightarrow H^*K(A, n-1).
$$
Define $K$ to be the Hopf cokernel of the morphism $i^*\colon H^*E
\rightarrow H^*F$ and $M$ to be the kernel of the morphism $\pi^*$
on primitive elements $PH^*K(A, n) \rightarrow PH^*E$.  The
cohomology suspension morphism $\sigma: H^*(K(A,n))\rightarrow
\Sigma H^*(K(A,n-1))$ restricted to $M$ defines a morphism
$M\rightarrow \Sigma K$. The adjoint $\Omega M \rightarrow K$
induces an isomorphism of $\A_p$-Hopf algebras $U \Omega M
\rightarrow K$ (compare with \cite[Proposition~5.7]{MR34:8406}),
where $U$ is Steenrod-Epstein's functor, left adjoint to the
forgetful functor $\K \rightarrow \U$.

Denote by $N$ the cokernel of $M \hookrightarrow PH^*K(A, n)$. The
above complex is then exact at $H^*E$,
\cite[Proposition~5.5]{MR34:8406}, and its homology at $H^*F$ is
isomorphic to $U\Omega_1 N$, where $\Omega_1$ is the first left
derived functor of $\Omega$.

Theorem~\ref{complement} now follows from the fact that $\Omega_1
N$ is a finitely generated unstable module. This simply reflects
the fact that the functor $\Phi$ (which is the ``doubling" functor
when $p=2$), \cite[1.7.2]{MR95d:55017} takes finitely generated
unstable modules to finitely generated ones.}
\end{remark}

We have thus shown that Eilenberg-MacLane spaces $K(A,n)$ satisfy
the strong fg closure property. But in fact, it can be easily
generalized to $p$-torsion Postnikov pieces.

\begin{proposition}
\label{Postnikovfiber}
Consider an $H$-fibration $F \xrightarrow{i} E \xrightarrow{\pi}
B$, where $B$ is an $p$-torsion $H$-Postnikov piece whose homotopy
groups are finite direct sums of cyclic groups and Pr\"ufer
groups. Then $H^*(E)$ is a finitely generated $\A_p$-algebra, if
and only if so is $H^*(F)$.
\end{proposition}

\begin{proof}
An induction on the number of homotopy groups of $B$ with
Lemma~\ref{extension} reduces to the case when $B$ is an
Eilenberg-Mac Lane space $K(A, n)$. We know from
Theorem~\ref{complement} that the statement holds in this case.
\end{proof}

Our first corollary has already been proved in
\cite{deconstructing}.

\begin{corollary}
\label{Postnikovpiece}
Let $F$ be an $H$-Postnikov piece of finite type. Then $H^*(F)$ is
finitely generated as unstable algebra.
\end{corollary}

\begin{proof}
The result is true for an Eilenberg-Mac Lane space $K(A, n)$ where
$A$ is an abelian group of finite type. The proof then follows by
induction on the number of homotopy groups.
\end{proof}

\section{Closure properties of $H$-fibrations}
\label{section closure}
The aim of this section is to extend the results of the preceding
section to arbitrary base spaces. We will prove that any $H$-space
$B$ such that $H^*(B)$ is a finitely generated algebra over $\A_p$
satisfies the fg closure property. We need here some input from the
theory of localization. Recall (cf. \cite{Dror}) that, given a
pointed connected space $A$, a space $X$ is $A$-\emph{local} if
the evaluation at the base point in $A$ induces a weak equivalence
of mapping spaces $\map(A, X) \simeq X$. When $X$ is an $H$-space,
it is sufficient to require that the pointed mapping space
$\map_*(A, X)$ be contractible.

Dror-Farjoun and Bousfield have constructed a localization functor
$P_{A}$ from spaces to spaces together with a natural
transformation $l: X \rightarrow P_{A}X$ which is an initial map
among those having an $A$-local space as target (see \cite{Dror}
and \cite{MR57:17648}). This functor is known as the
$A$-nullification. It preserves $H$-space structures since it
commutes with finite products. Moreover, when $X$ is an $H$-space,
the map $l$ is an $H$-map and its fiber is an $H$-space.

Recall that, for any elementary abelian group $V$, tensoring with
$H^*V$ has a left adjoint, Lannes' $T$-functor $T_V$,
\cite{MR93j:55019}. When $V = \Z/p$, the notation $T$ is usually
used instead of $T_{\Z/p}$ and the reduced $T$-functor is left
adjoint to tensoring with the reduced cohomology of $\Z/p$. This
allows to characterize the Krull filtration of the category $\U$
of unstable modules as follows: $M \in \U_n$ if and only if
$\overline{T} ^{n+1}M=0$, \cite[Theorem~6.2.4]{MR95d:55017}.

\begin{lemma}
\label{nullification}
Let $X$ be an $H$-space such that $T_V H^*(X)$ is of finite type
for any elementary abelian $p$-group $V$. Then $H^*(P_{B\Z/p} X)$
is finite if and only if, for any $n$, $H^*(P_{\Sigma^n B\Z/p} X)$
is a finitely generated $\A_p$-algebra.
\end{lemma}

\begin{proof}
By \cite{B2}, for any $n$ there are fibrations $P_{\Sigma^n B\Z/p}
X \rightarrow P_{\Sigma^{n-1} B\Z/p} X\rightarrow K(A_n,n+1)$
where $A_n$ is a $p$-torsion abelian group, which is a finite
direct sum of copies of cyclic groups $\Z/p^r$ and Pr\"ufer groups
$\Z_{p^\infty}$ (the technical hypothesis on the $T$ functor allows
to apply \cite[Theorem~5.4]{deconstructing}). In this
situation, we can apply Theorem \ref{complement} to show that
$H^*(P_{\Sigma^n B\Z/p} X )$ is a finitely generated algebra if
and only if $H^*(P_{\Sigma^{n-1} B\Z/p} X )$ is so. The statement
follows by induction since $H^*(P_{B\Z/p} X)$ is always locally
finite.
\end{proof}

\begin{theorem}
\label{generalclosure}
Consider an $H$-fibration $F \xrightarrow{i} E \xrightarrow{\pi}
B$. If $H^*(F)$ and $H^*(B)$ are finitely generated
$\A_p$-algebras, then so is $H^*(E)$.
\end{theorem}

\begin{proof}
Since both $H^*(F)$ and $H^*(B)$ are finitely generated
$\A_p$-algebras, the modules of indecomposable elements $QH^*(F)$
and $QH^*(B)$ are finitely generated $\A_p$-modules. Therefore,
\cite[Lemma 7.1]{deconstructing}, they belong to some stage
$\U_{n-1}$ of the Krull filtration. By \cite[Theorem
5.3]{deconstructing}, we know that both $F$ and $B$ are
$\Sigma^{n} B{\Z/p}$-local spaces. Since nullification preserves
fibrations whose base space is local (see
\cite[Corollary~3.D.3]{Dror}), it follows that $E$ is also
$\Sigma^n B{\Z}/p$-local.

Let us consider the fibration $\bar P_{B\Z/p} B \rightarrow B
\rightarrow P_{B\Z/p} B$. We know from Lemma~\ref{nullification}
that the nullification $P_{B\Z/p} B$ has finite mod $p$
cohomology. The homotopy fiber $\bar P_{B\Z/p} B$ is an
$H$-Postnikov piece whose homotopy groups are finite direct sums
of cyclic groups $\Z/p^r$ and Pr\"ufer groups $\Z_{p^\infty}$ by
\cite[Theorem~5.4]{deconstructing}. Lemma~\ref{extension} and
Proposition~\ref{Postnikovfiber} show then that it is enough to
prove the theorem when $H^*(B)$ is finite.

In that case, $B$ is a $B\Z/p$-local space, and we have thus a diagram of
fibrations
\[
\diagram F \rto \dto & E \dto \rto & B \ar@{=}[d] \cr P_{B\Z/p} F
\rto & P_{B\Z/p} E \rto & B
\enddiagram
\]
The mod $p$ cohomology $H^*(P_{B\Z/p} F)$ is finite and hence so
is $H^*(P_{B\Z/p} E)$ by an easy Serre spectral sequence argument.
Finally, we can apply Lemma \ref{nullification} to conclude that
$H^*(E)\cong H^*(P_{\Sigma^n B\Z/p}E)$ is a finitely generated
$\A_p$-algebra.
\end{proof}

\begin{corollary}
\label{covers}
Consider an $H$-space $X$ with finite mod $p$ cohomology. Then the
mod $p$ cohomology of its $n$-connected cover $X\langle n \rangle$
is a finitely generated $\A_p$-algebra.
\end{corollary}

\begin{proof}
Consider the $H$-fibration $\Omega(X[n]) \rightarrow X\langle n
\rangle \rightarrow X$. The fiber is an $H$-Postnikov piece of
finite type and the cohomology of the base is finite. The result
follows.
\end{proof}

This can be seen as the mirror result of \cite{spinoff}, where we
proved that any $H$-space with finitely generated cohomology as
algebra over the Steenrod algebra is an $n$-connected cover of an
$H$-space with finite mod $p$ cohomology, up to a finite number of
homotopy groups.

\begin{remark}
\label{Krull}
{\rm We have already seen in the proof of
Theorem~\ref{generalclosure} that the module of indecomposable
elements of a finitely generated cohomology, as algebra over the
Steenrod algebra, belongs to some stage of the Krull filtration.
In fact, for any $H$-space $X$ with finite mod $p$ cohomology,
$QH^* X\langle n \rangle$ belongs to $\U_{n-2}$ by \cite[Theorem
5.3]{deconstructing}, because $\Omega^{n-1} (X\langle n \rangle)$
is $B\Z/p$-local.


}
\end{remark}

\begin{remark}
\label{iff}
{\rm Theorem~\ref{generalclosure} cannot be improved to an ``if
and only if" statement. Consider for example the path-fibration
for the $3$-dimensional sphere $\Omega S^3 \rightarrow PS^3
\rightarrow S^3$. It is well-known that $H^*(\Omega S^3)$ is a
divided power algebra, which is not finitely generated over
$\A_p$.}
\end{remark}

We conclude with a characterization of the $H$-spaces which
satisfy the strong fg closure property.

\begin{proposition}
\label{strongorweak}
Let $X$ be an $H$-space $X$ which satisfies the strong fg closure
property. Then $X$ is, up to $p$-completion, a $p$-torsion
Postnikov piece.
\end{proposition}

\begin{proof}
If $X$ satisfies the strong fg closure property, then $H^*(\Omega
X)$ is a finitely generated $\A_p$-algebra. But in this case, by
\cite[Corollary~7.4]{deconstructing}, $\Omega X$ is, up to
$p$-completion, a $p$-torsion Postnikov piece.
\end{proof}





\bigskip
{\small
\begin{center}
\begin{minipage}[t]{6 cm}
Nat\`{a}lia Castellana and J\'er\^{o}me Scherer\\
Departament de Matem\`atiques,\\
Universitat Aut\`onoma de Barcelona,\\
E-08193 Bellaterra, Spain\\
\textit{E-mail:}\texttt{\,natalia@mat.uab.es}, \\
\phantom{\textit{E-mail:}}\texttt{\,\,jscherer@mat.uab.es}
\end{minipage}
\begin{minipage}[t]{5 cm}
Juan A. Crespo\\
Departamento de Econom\'\i a,\\
Universidad Carlos III de Madrid,\\
E-28903 Getafe, Spain\\
\textit{E-mail:}\texttt{\,jacrespo@eco.uc3m.es},
\end{minipage}
\end{center}}

\end{document}